\documentclass[11pt,leqno]{amsart}
\usepackage{t1enc}
\usepackage[latin1]{inputenc}
\usepackage[english]{babel}

\usepackage{amsthm, amssymb,color,wasysym,multirow, mathtools}
\usepackage{yfonts}

\usepackage{bbm}
\usepackage{bm}
\usepackage{mathrsfs}
\usepackage[normalem]{ulem}
\usepackage{graphicx}

 \usepackage{hyperref}
 \hypersetup{
     colorlinks=true,
     linkcolor=blue,
     filecolor=magenta,
citecolor = blue,
     urlcolor=blue,
 }

\usepackage[cmtip,all]{xy}
\newcommand{\longsquiggly}{\xymatrix{{}\ar@{~>}[r]&{}}}

\newtheorem{conjecture}{Conjecture}

\newcommand{\seteq}{\coloneqq}

\newcommand{\be}[0]{\begin{equation}}
\newcommand{\ee}[0]{\end{equation}}

\newcommand{\blue}{\color{blue}}

\newcommand{\bit}{\begin{itemize}}
\newcommand{\eit}{\end{itemize}}
\newcommand{\ber}{\begin{red}}
\newcommand{\er}{\end{red}}
\newcommand{\beb}{\begin{blue}}
\newcommand{\eb}{\end{blue}}
\newcommand{\ed}{\end{document}}

\newcommand{\CD}{{\mathcal D}}
\newcommand{\CP}{{\mathcal P}}
\newcommand{\CS}{{\mathcal S}}

\numberwithin{equation}{section}

\theoremstyle{plain}

\theoremstyle{definition}

\begin{document}

\title[Mathematical Data Science]{Mathematical Data Science}

\author[M. R. Douglas and K.-H. Lee]{Michael R. Douglas and Kyu-Hwan Lee}


\maketitle

\section{Introduction}

Can machine learning help discover new mathematical structures? 
In this article we discuss an approach to doing this which one can call {\em mathematical data science}.  In this paradigm, one studies mathematical objects collectively rather than individually, by creating datasets and doing machine learning experiments and interpretations.  Broadly speaking, the field of data science is concerned with assembling, curating and analyzing large datasets, and developing methods which enable its users to not just answer predetermined questions about the data but to explore it, make simple descriptions and pictures, and arrive at novel insights.  This certainly sounds promising as a tool for mathematical discovery! 

Mathematical data science is not new and has historically led to very important results.
A famous example is the work of Birch and Swinnerton-Dyer leading to their conjecture \cite{BSD},
based on computer generation of elliptic curves and linear regression analysis of the resulting data.
However, the field really started to take off with the deep learning revolution and with the easy access to ML models
provided by platforms such as {PyTorch} and {TensorFlow}, and built into computer algebra systems such as {Mathematica, Magma} and {SageMath}.
Some early influential works include \cite{HE2017564,Carifio2017,Ruehle2017,Krefl2017}.\footnote{
Many of these works were done by physicists, who have the advantage over mathematicians
that they can more easily publish conjectures. }  
Of particular note is \cite{davies2021advancing} which brought more ML
sophistication (such as attribution analysis) and, importantly, was done in collaboration with expert mathematicians who could
guide the search for interesting conjectures and then prove them.  After all, even if a work makes intriguing claims about 
mathematical objects, it is mathematics only to the extent that the claims (or others inspired by it) can eventually be rigorously proven.

With the continuing increase of computer power and software sophistication, the development of mathematical databases, 
and the wealth of ML methods now available,
the prospects for ML approaches look very bright.  On the other hand, there are also difficulties which are particularly salient for mathematics: the challenge
of interpreting ``black box'' ML models, and especially the need to turn evidence into precisely formulated conjectures and rigorous proofs.

After an overview,
we will present two case studies from the work of the second author: murmurations in number theory and loadings of partitions related to Kronecker coefficients in representation theory and combinatorics. 

\subsection{The mathematical data science (MDS) paradigm}
\label{ss:mds}

This 
can be summarized as follows.
\begin{enumerate}
 \item Consider mathematical objects collectively to generate datasets. 
 \item Apply ML tools to find structure (usually statistical) in the datasets.
 \item Interpret these results to understand the mathematical objects better. 
 \item Find new precisely defined structures, conjectures and theorems.
\end{enumerate}
The first comment to make is that (to date) all of these steps are carried out by human mathematicians,
employing the computer as a tool.  They all involve choices and intuition which we do
not know how to automate.  Thus we are not in danger of being replaced by the computer (yet).

Let us now make more detailed comments about each step.
First, the data.  This is (of course) synthetic data, generated by a computer following a mathematically precise definition.
Typically a dataset is a list of objects together with values of associated invariants (the ``features'' in ML parlance).
This requires making choices: when are two objects isomorphic; what invariants to include.
In some fields (for example number theory and knot theory) there is a clear sense for what objects and invariants are interesting to list
and what isomorphisms to impose, efficient algorithms have been developed to do this, and datasets have been
computed and made publicly available (e.g., \cite{lmfdb,regina}).
However this is by no means required and the first step for many projects is to generate a new dataset.  
An example is the dataset of Kronecker coefficients in \S \ref{s:kron}.

A very important point which is often passed over too quickly is that the ML dataset is usually {\bf not} the actual set one wants
to understand.  Rather, it is a sample from some larger ``population'' of objects which is the true subject of study.  
For example, we might intend to study knots, but of course there are an infinite number of distinct knots and no finite database
can contain them all.  One must choose a finite subset and hope that it suffices for the questions under study.  
If we are studying statistical properties of the set (its distribution or that of some invariant, averages, {\it etc.}), 
we need a representative subset, one for which the properties being measured are a good estimate of the properties
of the entire population.


Often one can define a finite subset by imposing a mathematically natural criterion, such as all objects up to a certain 
complexity or size.
For example, most tables of knots list all knots up to a specified maximal crossing number.
This use of a size parameter (such as crossing number) is common and will be used in the case studies below,
but it is not the only possibility.  Another possibility is to define
a probability distribution on the full set, and create the dataset by sampling from this distribution.  
Why would one do this?  One practical reason is just to reduce the size of the dataset to what the computer at hand can work with.
But there are other reasons.  One is that such a sampled dataset can actually be more representative (in a statistical sense)
than simply taking all objects up to a specified size. Another is that
often the challenge in uniquely enumerating objects (one representative
of each isomorphism class) is not to generate the objects or compute the invariants, it is to identify the isomorphism classes
(this is the case for knots \cite{burton:LIPIcs.SoCG.2020.25}).  A sampling algorithm which generates multiple (and varying numbers of)
representatives of every class can be much simpler (see for example \cite{even-zohar_models_2017} for knots).  One still needs to know something
about the probabilities, but they need not be uniform.  Of course, using sampling will generate random datasets,
but this randomness can be defined mathematically and controlled.



The essential point of a mathematical dataset is not so much that it was computed and is available for download and use, convenient as this may be.
Rather it is that it has a simple and mathematically precise definition, so that the results it leads to can be stated precisely.
We offer the term ``platonic dataset'' for such datasets;
let us define one as:
\bit
\item A set $\CS$ of mathematical objects.
\item A function $F$ on $\CS$ into some space of invariants $\CD$, say ${\mathbb R}^d$.
\item A way to pick out finite subsets $\CS_i\subseteq\CS$ -- this could be
filtration by size, or sampling from a probability distribution $\rho$ on $\CS$.
\eit
Given this definition, we can apply $F$ to each $\CS_i$ to get a distribution $\CD_i$ on $\CD$, usually a point cloud 
(a subset consisting of a finite number of points) in practice.


Next, let us discuss the ML tools, at a very abstract level.  
The upshot of the previous discussion is that we have some datasets which, from the point of view of ML,
are point clouds $\CD_i\subset \CD$.  We can apply to these two standard paradigms of ML, supervised learning and unsupervised (or self-supervised) learning. 
In self-supervised learning, we model the probability distribution $\rho_i$ of the population from which $\CD_i$ was sampled.
In supervised learning, we split $\CD$ into ``input'' and ``output'' variables, $\CD \cong \CD_{in}\times \CD_{out}$, and find
a function $y_i:\CD_{in}\rightarrow \CD_{out}$ which models the subset of $\CD$ on which $\CD_i$ has support.  In both cases by ``model''
we mean, out of a parameterized class of distributions $\rho_\theta$ or functions $y_\theta$, we find one which has a small error
measured by some objective (or ``loss'') function: say cross-entropy (for probability), or mean squared error (for function fitting).
Modern ML provides very general and flexible function classes (feed-forward networks, transformers, {\it etc}.) and can do these operations
on datasets with billions of items, even given only a laptop with a GPU.

In addition there are numerous auxiliary tools of data science -- topological data analysis, visualization 
tools and the like -- and a good deal of statistical and other theory and intuition we can call on.  
This includes many ways to deal with the problem mentioned above, that ML datasets are often only
subsets of the complete set of objects of interest, and might or might not be representative of the properties we want to study.
In ML terms one often speaks of a model working ``in distribution'' (so, on data with the same statistics as the
dataset used for training the model) but not ``out of distribution'' (on the actual population).
This is very common, leading to the universal practice of testing models by evaluating the error on a separate dataset
of objects ``held out'' for this purpose.  These concepts lie behind the important practical wisdom and advice from the ML experts, which can be of immense  value in carrying out a mathematical data science project.

While this brief summary of modern ML and data science may seem too reductive, and certainly it leaves out a great deal,
nevertheless we present it as a reasonable first picture of what ML brings to the table.  What good is it for a mathematician?
First of all, it leads one to think about one's problems in a different way.  As we will see in \S \ref{s:mur}, just by looking at
the much studied subject of $L$-functions of elliptic curves through this lens, interesting discoveries were made.  Next, although
the patterns which one might find using these techniques are often rather simple, such as inequalities bounding the dataset,
they can still be important.  We will give examples in \S \ref{s:kron}.

Steps 3 and 4 in our outline, to interpret the results in terms of structure, 
are of course where the discoveries are made.  These are best discussed using case studies.

\section{Murmurations} 
\label{s:mur}

Here is a picture of murmuration in nature. In mathematics, murmuration was first discovered in a dataset of elliptic curves \cite{HLOP}.

\begin{center}
\includegraphics[scale=0.4] {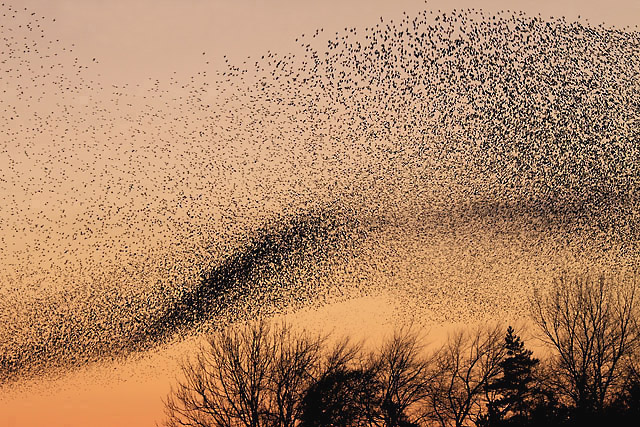}
\end{center}
{\tiny (Photo by Walter Baxter / A murmuration of starlings at Gretna / CC BY-SA 2.0)}

\subsection{Elliptic Curves}

An elliptic curve over $\mathbb Q$ can be defined by a cubic equation $y^2 = x^3 + a x + b$ with integer coefficients $a$ and $b$;
pairs $(x,y)\in \mathbb Q\times\mathbb Q$ satisfying the equation are called rational points.
Consider 
\be y^2=x^3+x+1 ,
\ee
it has label 496.a1 in the LMFDB database of elliptic curves \cite{lmfdb}.
There are infinitely many distinct elliptic curves; one could define a finite subset by $|a|\le N,|b| \le N$.
The LMFDB uses a related but more sophisticated measure of the complexity of a curve called the {\em conductor}; 
it includes all curves of conductor less than $500{,}000$, and many more.

The Mordell--Weil theorem tells us that the rational points of an elliptic curve form a finitely generated abelian group. This has a torsion part and a free part, and the rank of the free part is called the {\em rank} of the elliptic curve. Elliptic curves have been studied for a long time, but the rank is still mysterious, with many questions still unanswered: in particular, the Birch and Swinnerton-Dyer conjecture which we will come to below.

Can we apply ML to classify elliptic curves according to their rank? The first step is to decide how to represent an elliptic curve to a machine. 
While the data $(a,b)$ suffices in principle, there are many other invariants which are nontrivial functions of this data, and it may be better to
include them as well.

A very successful way to study arithmetic objects is through their $L$-functions.
Let $E$ be an elliptic curve over $\mathbb Q$. For each prime $p$, we can reduce the curve modulo $p$ and count the number of rational points. We denote this number by $n_p$, and define $a_p= p+1 -n_p$. These $a_p$'s are used to define the $L$-function of the elliptic curve $E$ as 
\be\label{eq:defLfn}
L(E,s) = \prod_{\substack{p\ \small\mbox{prime}\\ p \nmid N_E}} \frac{1}{1-a_p(E) p^{-s} + p^{1-2s}}\prod_{\substack{p\ \small\mbox{prime}\\ p \mid N_E}} \frac{1}{1-a_p(E) p^{-s}},
\ee
where $N_E$ is the conductor of $E$.
While this only converges for $\mbox{Re}(s)>3/2$, 
it turns out to be a meromorphic function of $s$, allowing it to be defined by analytic continuation elsewhere.
According to the Birch and Swinnerton-Dyer conjecture, the order of its vanishing at $s=1$ is the rank of the curve.
The conjecture also predicts the coefficient of the expansion around $s=1$ as a product of invariants of the curve.
By taking logarithms one gets a linear relation, and the original evidence for it came from linear regression
on a dataset of curves.

Rather than attempt to make direct contact with this, we fix a set of $k$ primes, compute the $a_p$'s for these, and use the sequence 
$\{a_{p_1},a_{p_2},\ldots,a_{p_k}\}$
to represent the elliptic curve.  We can think of it as a ``picture'' of the curve.  Since we want to predict the rank, we adjoin the rank 
to our dataset as well.
Thus, each curve corresponds to a point in a space of invariants $\CD$ which is the product of $\CD_{in} = {\mathbb R}^k$ and $\CD_{out}$ representing the rank.
For reasons which are standard in ML, we take $\CD_{out}\cong {\mathbb R}^n$ and represent a curve of rank $0\le r<n$ as the $r$'th basis vector 
$(0^r,\ 1,\ 0^{n-r-1})$.%
\footnote{One argument for doing this is that the model will actually be predicting the vector of probabilities $(P_0,P_1,\ldots)$ that a curve
with input $\{a_p\}\in \CD_{in}$ has rank $r$.}

In collaboration with He, Oliver and Pozdnyakov, the second author created a dataset consisting of 16,000 elliptic curves of rank 0 and 16,000 curves of rank 1, with conductors ranging from 1 to 10,000. They chose the number of primes $k= 300$ and computed the dataset as above.  Since the ML problem is to partition the items into two classes (rank 0 and 1), they applied the standard technique for doing this, called logistic regression.

The results showed an impressive accuracy of nearly 99\%. When they reduced the number of $a_p$'s,' the accuracy dropped slightly, but it remained very high. In other experiments of classifying elliptic curves with ranks 0, 1, and 2, the outcomes were similar. In particular, distinguishing between rank 1 and rank 2 curves was the easiest task for the ML model.  All this was a satisfying demonstration of the effectiveness of machine learning. 

\subsection{Principal Component Analysis}

In addition to supervised learning (training using the predicted quantity), one might ask whether there is some
{\it a priori} structure in the $a_p$'s which determines the rank -- a question for unsupervised learning.
The logistic regression which was found to be effective 
is a variation on linear regression, showing that the problem is amenable to a linear analysis.

The simplest but still most important tool for unsupervised learning is linear: it is principal component analysis (PCA).
Given a point cloud $\{x_i\}\subset \CD\cong{\mathbb R}^d$, we define the empirical covariance matrix $(XX) \in \CD\otimes\CD$ as
\be
(XX) \equiv \frac{1}{N}\sum_{i=1}^N x_i \otimes x_i .
\ee
PCA is then simply a reduction of this matrix along its principal axes,
\be
(XX) = R \, \Lambda \, R^t
\ee
where $R$ is an orthogonal transformation and $\Lambda$ is a diagonal linear transformation (the eigenvalues of $(XX)$).
Often one regards the largest eigenvalues as the most important and truncates to a finite number $n$ of these -- below we
will just keep the largest $n=2$ eigenvalues.  The $m$'th principal component of a datapoint $x_i$ (called PC$m$) is then 
its projection onto the corresponding eigenvector, $x_i \cdot R \cdot e_m$.  In the language of \S \ref{ss:mds}, we
are modeling the data distribution as a normal (Gaussian) distribution supported on an $n$-dimensional subspace, and
looking at the projection of the data on this subspace.

Here is the result of PCA on the dataset of rank 0 and rank 1 elliptic curves. There is a clear separation along the PC1 axis.

\begin{center}
     \includegraphics[scale=0.17]{"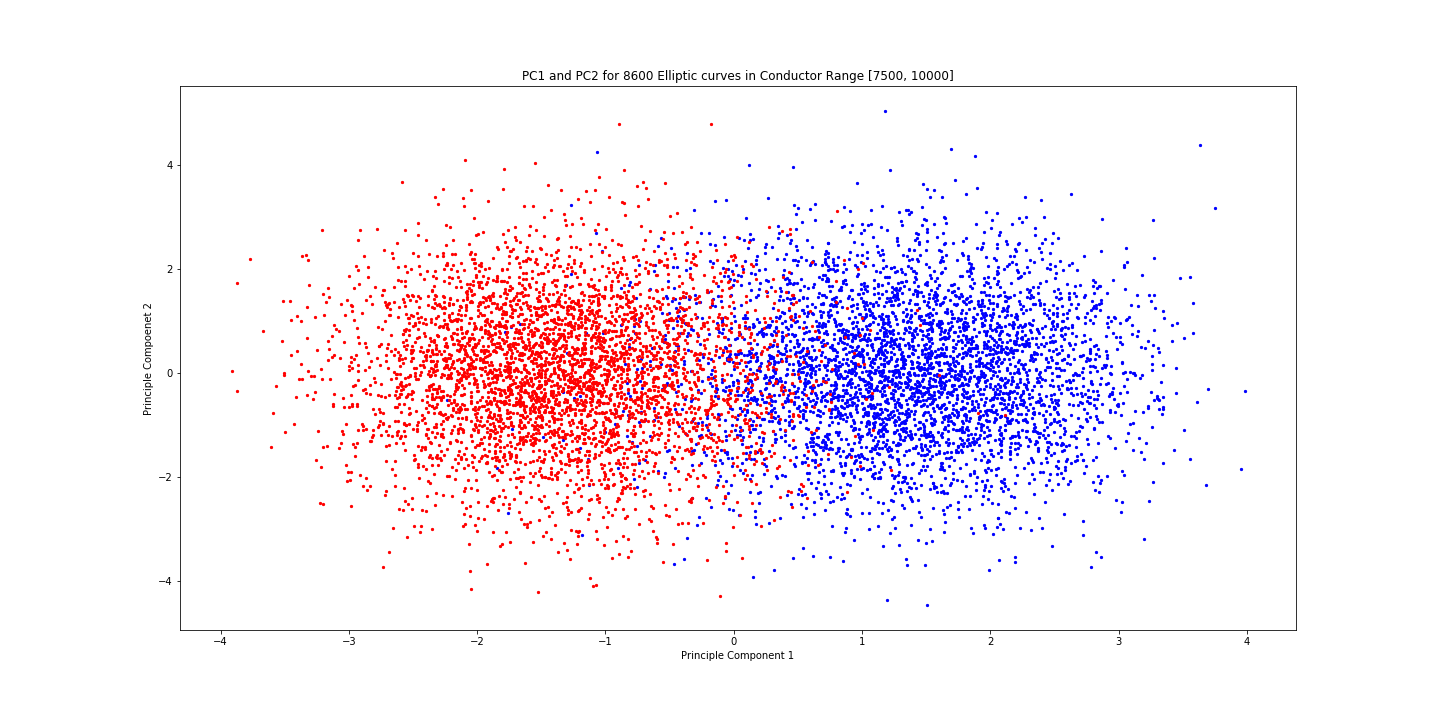"}
\end{center}
 {\scriptsize \sf A plot of PC1 against PC2 for the balanced dataset of $8,600$ randomly chosen elliptic curves with rank $r_E\in\{0,1\}$ and conductor $N_E\in[7500,10000]$. The blue (resp. red) points correspond to curves with rank $0$ (resp. $1$).}

\vskip 0.1 cm

Here is the result of PCA on the dataset of rank 0, 1 and 2 elliptic curves. The green cluster corresponds to the rank 2 curves, and it is clearly separated from the other clusters, again along the PC1 axis.

\begin{center}
     \includegraphics[scale=0.19]{"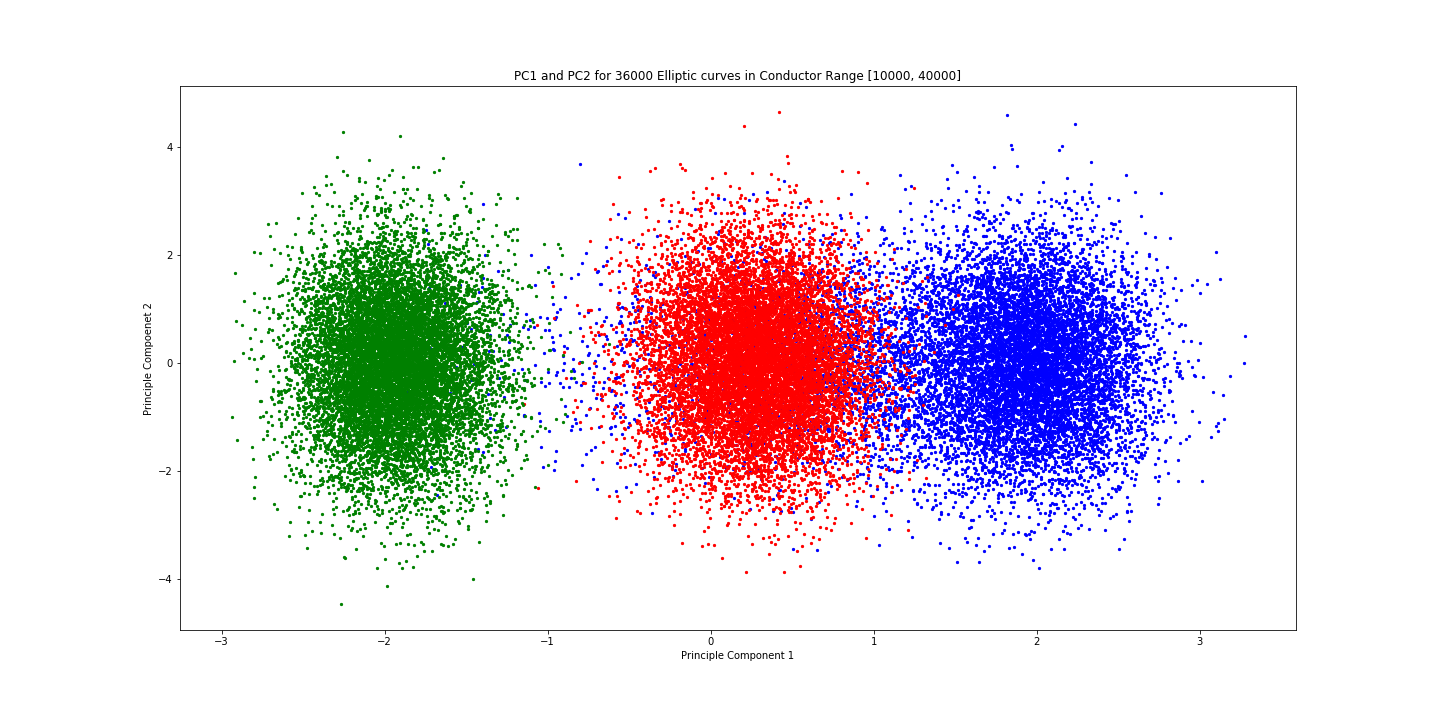"}
\end{center}
{\scriptsize \sf A plot of PC1 against PC2 for elliptic curves in the balanced dataset of $36,000$ randomly chosen elliptic curves with rank $r_E\in\{0,1,2\}$ and conductor $N_E\in[10000,40000]$. The blue (resp. red, green) points correspond to curves with rank $0$ (resp. $1$, $2$). }

\vskip 0.1 cm

Here is a picture which shows the components of PC1 (the top eigenvector of $(XX)$) for the rank 0 and rank 1 case. 
We took curves with conductor $N_E\in[1000,2000]$.
The $x$-axis shows the prime count. What is interesting is that the $y$-values are oscillating.

\begin{center}
     \includegraphics[scale=0.4]{"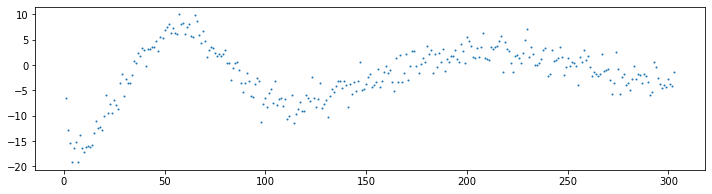"}
\end{center}

\subsection{Murmurations} 

This oscillation was interesting.  What caused it?  Was it something about PCA, or was it a simpler property of our data?  
To address this, we fix a prime $p$ and compute the average of $a_p$ over a set of elliptic curves for each rank with conductors in a fixed interval. This gives us a function where $n$ is the prime count and $r$ is the rank and the interval runs from $N_1$ to $N_2$,
\begin{equation*}
f_r(n) =\frac{1}{\#\mathcal{E}_r[N_1,N_2]}\sum_{E\in\mathcal{E}_r[N_1,N_2]}a_{p_n}(E),
\end{equation*}
where $N_1 < N_2\in\mathbb{Z}_{>0}$, and $\mathcal{E}_r[N_1,N_2]$ is the set of elliptic curves over $\mathbb{Q}$ with rank $r$ and conductor in range $[N_1,N_2]$.
We are literally averaging $a_p$'s. Finally, we plot $f_0(n)$ in blue and $f_1(n)$ in red 
 for $1 \le n \le 1000$ and $[N_1,N_2]=[7500,10000]$. 

\begin{center}
     \includegraphics[scale=0.2]{"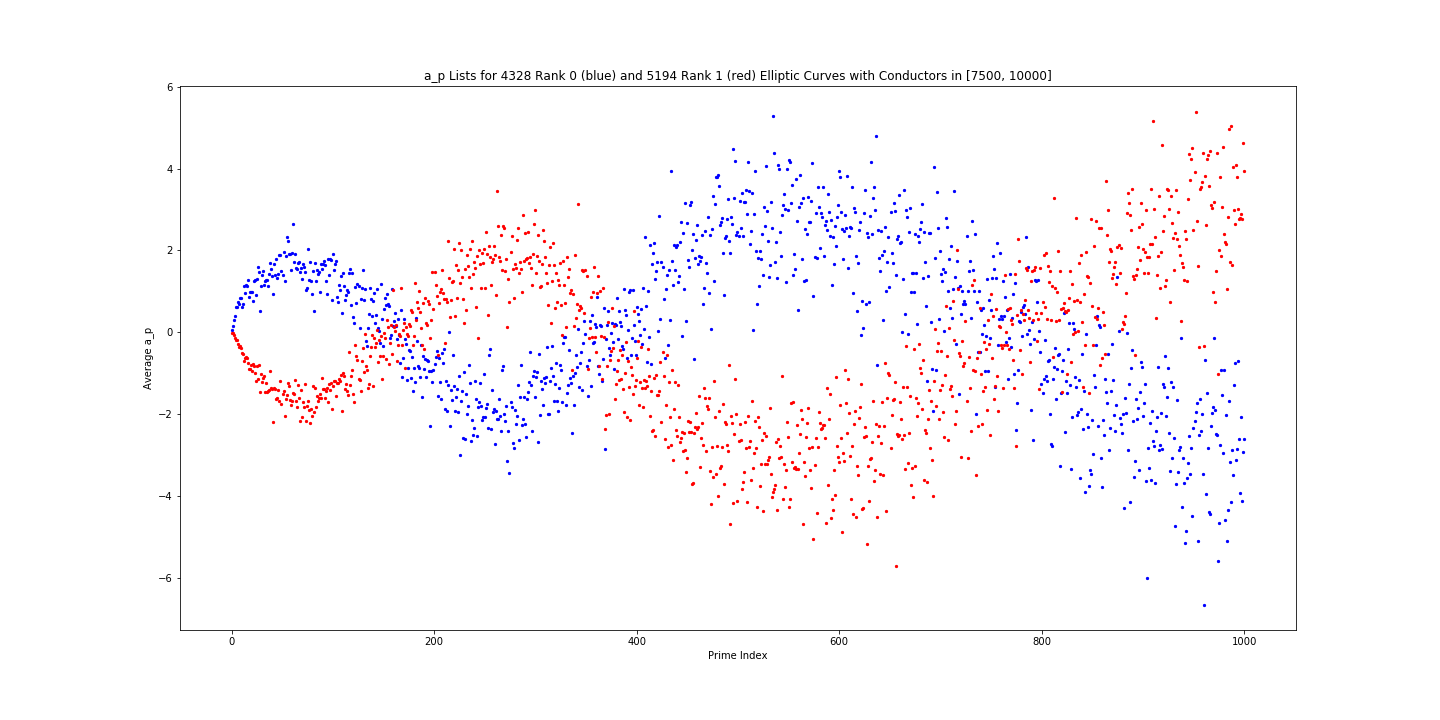"}

\end{center}
 
The first pattern is not surprising. We know that the $a_p$'s tend to be smaller when the rank is bigger. However, no one expected the averages to switch positions, and the oscillation came as a complete surprise.

Here is the plot for rank 0 and rank 2 curves. We can see that the oscillation patterns are similar. It suggests that the parity of the rank is important. 

\begin{figure}[h!!!]
\begin{center}
     \includegraphics[scale=0.23]{"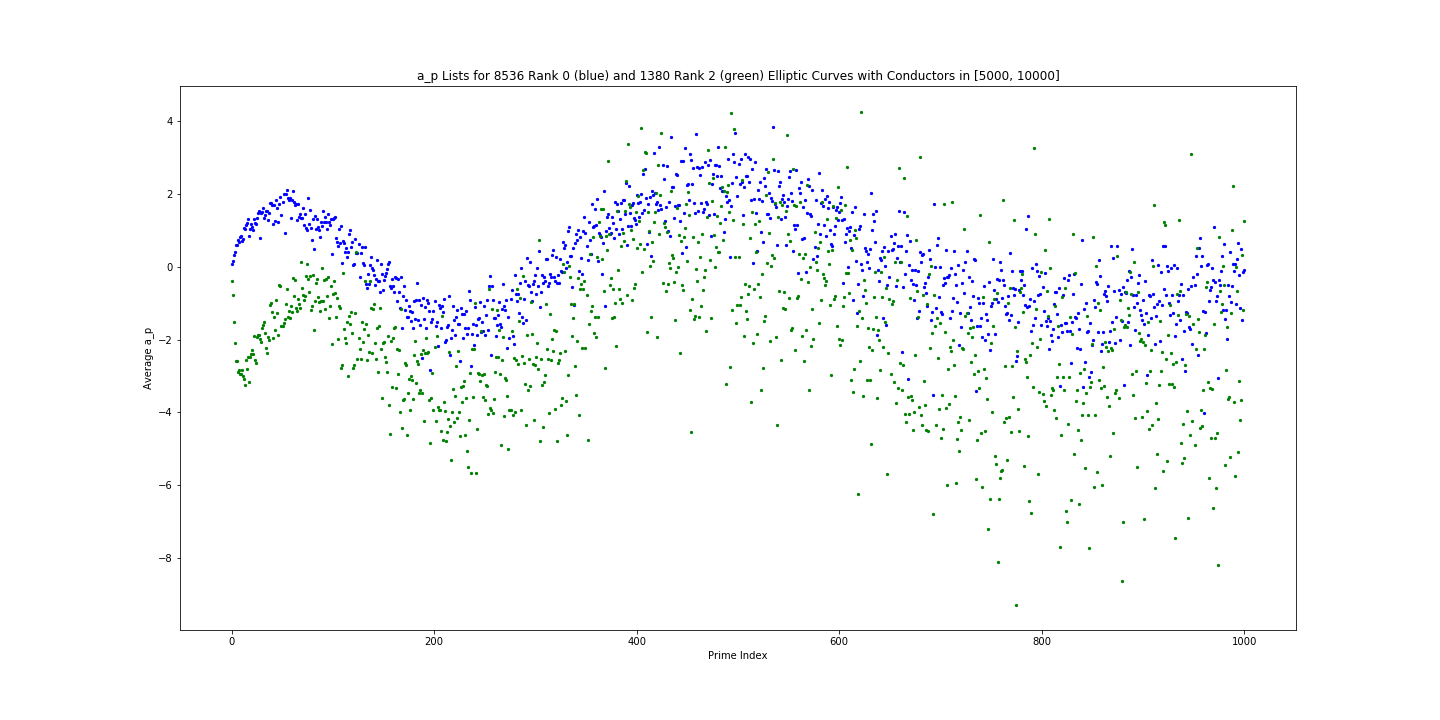"}
\end{center}
{\scriptsize \sf Plots of the functions $f_0(n)$ (blue) and $f_2(n)$ (green) for $1 \le n \le 1000$ and $[N_1,N_2]=[5000,10000]$.} 
\end{figure}

Next we ask: How does this oscillation change as the conductor range changes? To explore this, we used geometric intervals from $2^n$ to $2^{n+1}$ for the conductors and produced plots according to the parity of the rank. That is, we considered these functions and sketched the corresponding plots. At this point Andrew Sutherland joined us, and we thank him for the following plots.

Here are the resulting pictures for $n=13,14,15,16,17$. They are generated from distinct datasets, but they all show the same pattern. As we can see, the murmuration phenomenon appears to be scale-invariant.

\begin{center}
\includegraphics[width=4.5 in]{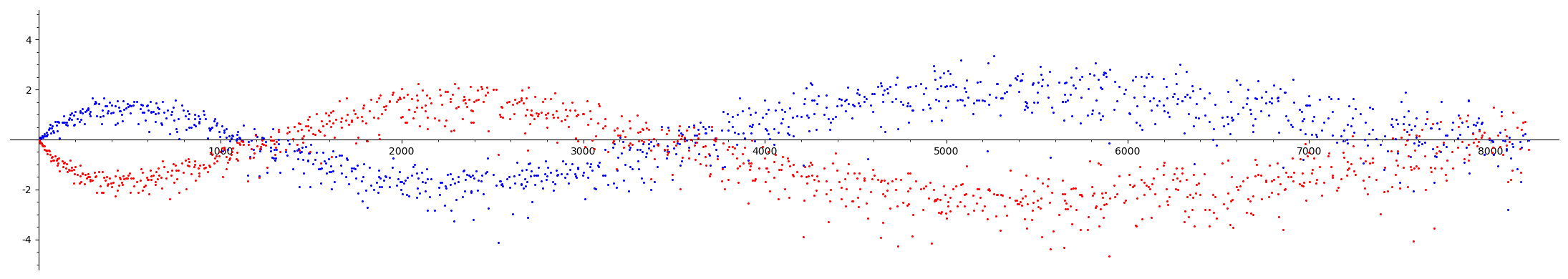}

\includegraphics[width=4.5 in]{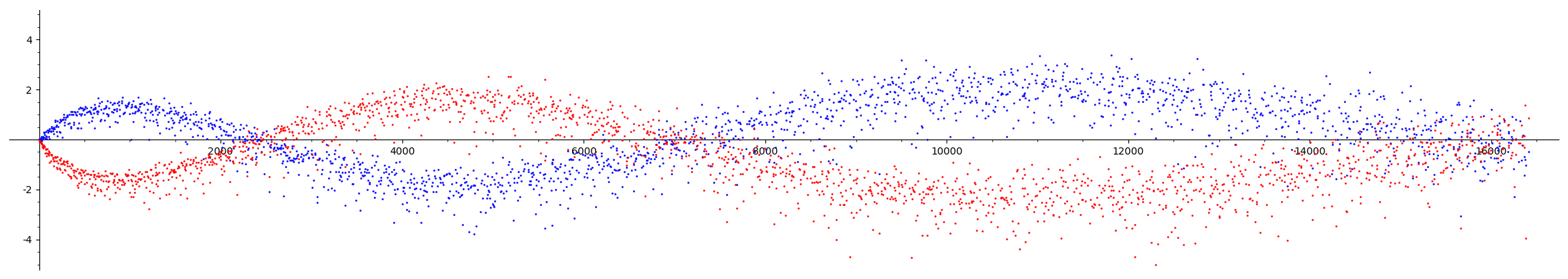}

\includegraphics[width=4.5 in]{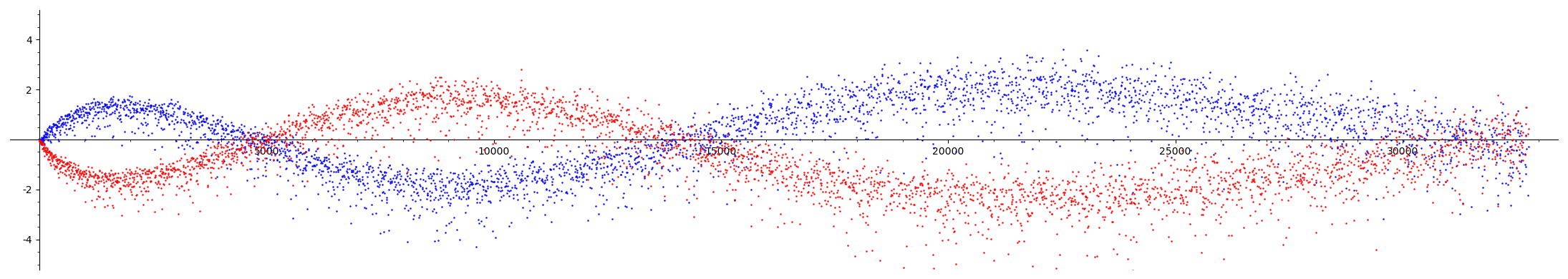}

\includegraphics[width=4.5 in]{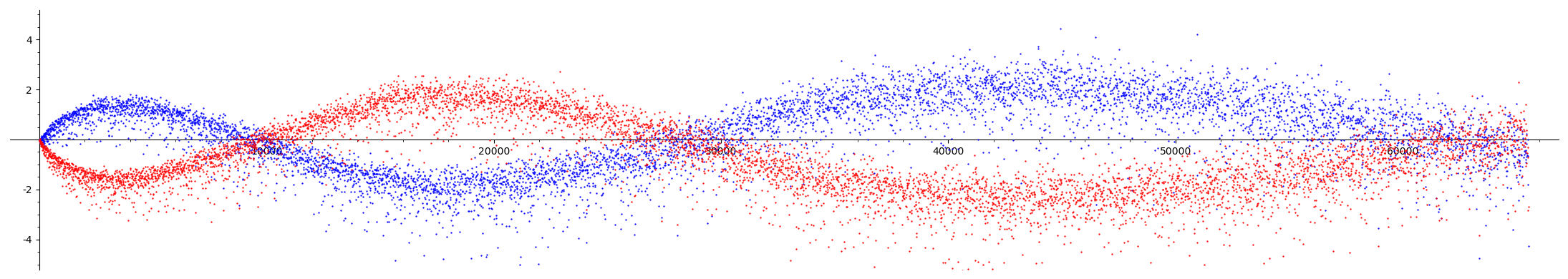}

\includegraphics[width=4.5 in]{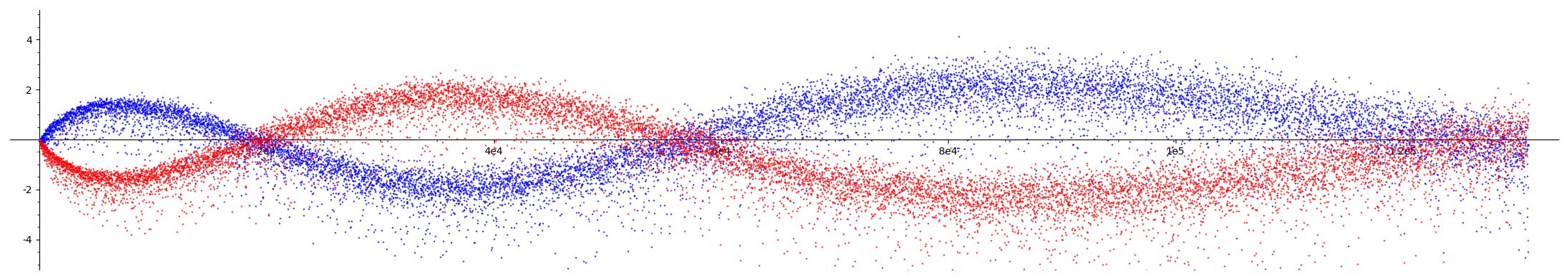}
\end{center}

The next question was: Does murmuration occur only for elliptic curves?  According to the modularity theorem, each elliptic curve over $\mathbb Q$ corresponds to a weight 2 modular form. However, there are modular forms of other weights, as well as higher genus curves. So, we can ask: Do murmurations appear for these objects too? The answer is yes. Andrew Sutherland conducted extensive computations and maintains a webpage dedicated to murmurations.\footnote{
 \url{https://math.mit.edu/~drew/murmurations.html} }
Here is a snapshot of the webpage.
\vskip 0.3 cm 
\includegraphics[width=4.8 in, height= 6 in]{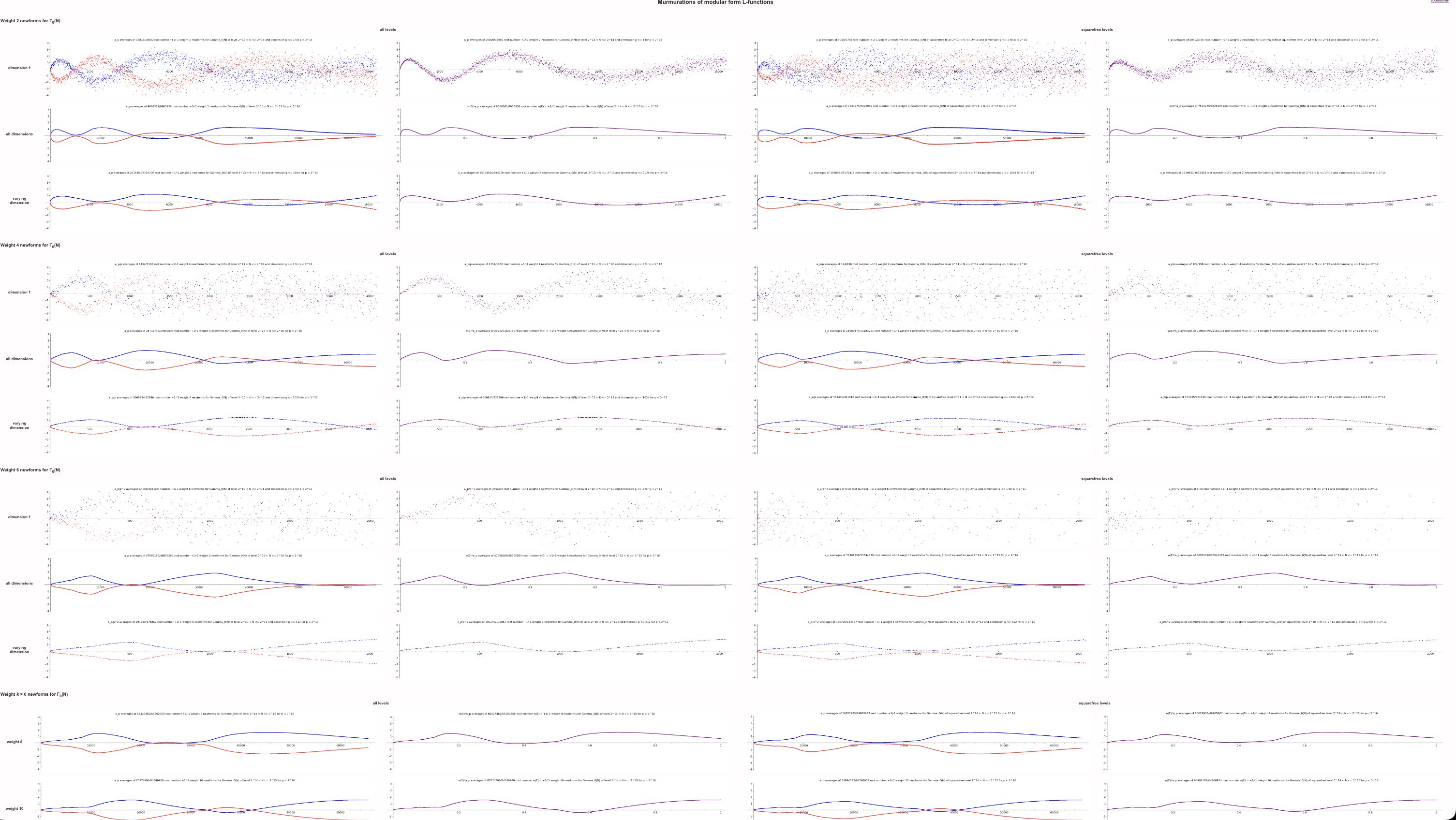}
\vskip 0.3 cm 

News of this discovery spread. In July 2023, there was a hot topics workshop on murmurations at ICERM, and in fall 2024 another workshop focused on murmurations at the Simons Center at Stony Brook University. The discovery was featured in Quanta magazine.
Number theorists began to look for explanations, rigorous proofs, and extensions of the basic idea.
Peter Sarnak from Princeton and IAS wrote notes on murmurations, and his student Nina Zubrilina proved the  murmurations for modular forms \cite{Z}.
Since then, many other papers on the topic have appeared on the arXiv, and research on murmurations is ongoing.  We see that machine learning experiments can lead to the discovery of new structures in number theory.

\section{Kronecker Coefficients}
\label{s:kron}

Next we will explore an example from representation theory and algebraic combinatorics. Following the same paradigm as for murmurations, we'll study Kronecker coefficients, as done in \cite{Kr}. Unlike murmuration, there is no canonical picture. We asked ChatGPT to create a picture related to Kronecker coefficients, and here is the outcome.  

\begin{center}
\includegraphics[scale=0.16] {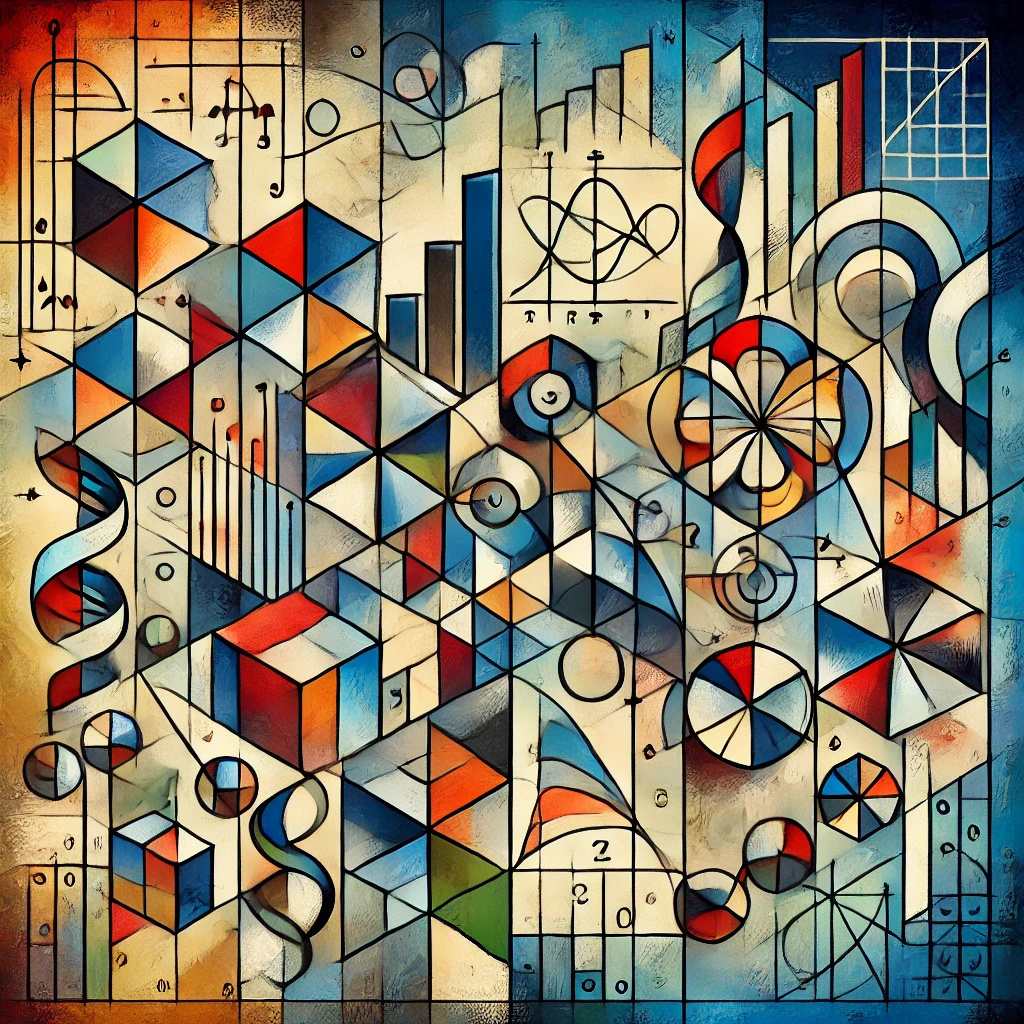}
\end{center}
{\tiny (generated by ChatGPT-4)}

\subsection{Kronecker coefficients}

Kronecker coefficients describe how a tensor product of two irreducible representations of the symmetric group decomposes into irreducible components. To fix the notations, let $\mathfrak S_n$ be the symmetric group of degree $n$. The irreducible representations of $\mathfrak S_n$ over $\mathbb C$ are parametrized by partitions $\lambda$. We denote them by $S_\lambda$, where $\lambda\in \CP(n)$ the set of partitions of $n$.  

Let's consider the tensor product of two irreducible representations, $S_\lambda$ and $S_\mu$. This product decomposes into a direct sum of irreducible representations,
\be
S_\lambda \otimes S_\mu = \bigoplus_{\nu \in \mathcal P(n)} g_{\lambda, \mu}^\nu S_\nu  \quad (g_{\lambda, \mu}^\nu \in \mathbb Z_{\ge 0}) .
\ee
As $\nu$ varies, the decomposition gives us the multiplicities $g_{\lambda, \mu}^\nu$. These multiplicities are called the Kronecker coefficients. They are fundamental structural constants in representation theory, but they are difficult to compute.

In the context of combinatorial representation theory, one can ask: Can we give a combinatorial description of the Kronecker coefficients? This question has remained open since 1938. It is known that deciding whether a Kronecker coefficient is non-zero is NP-hard \cite{ikenmeyer_vanishing_2017}.
And as in the case of the rank of an elliptic curve, we can ask: Can we use ML to determine whether $g_{\lambda,\mu}^\nu$ is nonzero or not?

\subsection{ML Kronecker coefficients}

Let $\CP(n)$ be the set of partitions of $n$, as before, and think of a partition of $n$ as an $n$-dimensional vector. For example, when $n$ is equal to 6, the partition $(5,1)$ becomes the 6-dimensional vector $(5,1,0,0,0,0)$. In this way, a triplet of partitions becomes a $3n$-dimensional vector. 
For example, when $n=14$, the dataset contains 42 dimensional vectors, and the total number of data points is approximately 2.46 million.  

Now for each $n$, we have an ML problem, where the dataset consists of these $3n$-dimensional vectors, and the target values are 0 and 1. 
Here 0 represents the case $g_{\lambda, \mu}^\nu$ is zero and 1 means it is nonzero. This setup makes it a binary classification problem. 
We could use $g_{\lambda, \mu}^\nu$ itself as the target value to frame the problem as multiclass classification. However, binary classification is already challenging, and in representation theory, it is often sufficient to determine whether $g_{\lambda, \mu}^\nu \neq 0$.

The machine learning results were very successful, with an accuracy as high as 98\%. For $n=12, 13, 14$, we applied nearest neighbors, convolutional neural networks and LightGBM, using various input formats. In particular, CNN$_2$ and CNN$_3$ use 2D and 3D arrays as inputs, respectively. All of these setups achieved high accuracy. 

\begin{center}
{\footnotesize \begin{tabular}{|c|c|c|c|c|c|}
\hline
\multirow{2}{*}{$n$} &\multirow{2}{*}{$\# \mathcal D$} & \multicolumn{4}{|c|}{Accuracy}\\
\cline{3-6} & & NearN  &CNN$_2$ & CNN$_3$&  LGBM \\ \hline 
$12$ & $ 126,900 \times 2$ & 0.9155 & 0.9529 & 0.9697 & 0.9714  \\
\hline
$13$ & $ 260,000 \times 2$ & 0.9318 & 0.9618 &0.9773& 0.9837  \\
\hline
$14$ & $ 600,000 \times 2$ & 0.9364 & 0.9635 & 0.9772 &  0.9845\\
\hline
\end{tabular}
}
\end{center}
\medskip
{\small \sf
The above table shows the accuracy of NearN, CNN$_2$, CNN$_3$ and LGBM classifiers when asked to distinguish $(\lambda, \mu, \nu)$ with $g_{\lambda, \mu}^\nu=0$ from those with $g_{\lambda, \mu}^\nu \neq 0$, where $\# \mathcal D$ is the number of sampled data points.} 

\subsection{Similitude matrix and $a$-loadings}

After completing these supervised learning experiments, we applied PCA and carefully analyzed the results.  This led us to explore different ways to embed the dataset into a higher dimensional space, thus introducing new concepts into the study of Kronecker coefficients. 

Let small $p(n)$ be the number of partitions of $n$, and $\mathsf P_n$ the $p(n)\times n$ matrix whose rows represent partitions. 
For example, 
\[{\scriptsize \mathsf P_6= \begin{bmatrix} 6&0&0&0&0&0 \\ 5& 1&0&0&0&0\\4& 2&0&0&0&0\\ 4& 1& 1&0&0&0\\ 3&3&0&0&0&0\\ 3& 2& 1&0&0&0\\3& 1& 1& 1&0&0\\2&2&2&0&0&0\\ 2& 2& 1& 1&0&0\\2& 1& 1& 1&1&0\\1&1& 1& 1& 1& 1\end{bmatrix} }.\] 
Each row represents a partition of 6, ordered in decreasing lexicographic order.  

We define the {\em similitude matrix} by the formula 
\be
\mathsf  Y_n \seteq \mathsf P_n \mathsf P_n^{\top}.
\ee
It is a symmetric integer matrix, with each entry indexed by a pair of partitions.  For example
\[  {\scriptsize \mathsf Y_6=\left [ \begin{array}{ccccccccccc}36&30&24&24&18&18&18&12&12&12& 6 \\ 30& 26& 22& 21& 18& 17& 16& 12& 12& 11& 6\\ 24& 22& 20& 18& 18& 16& 14& 12& 12& 10& 6\\ 24& 21& 18& 18& 15& 15& 14& 12& 11& 10& 6\\ 18& 18& 18& 15& 18& 15& 12& 12& 12&  9& 6\\ 18& 17& 16& 15& 15& 14& 12& 12& 11&  9& 6\\ 18& 16& 14& 14& 12& 12& 12& 10& 10&  9& 6\\ 12& 12& 12& 12& 12& 12& 10& 12& 10&  8& 6\\ 12& 12& 12& 11& 12& 11& 10& 10& 10&  8& 6\\ 12& 11& 10& 10&  9&  9&  9&  8&  8&  8& 6\\ 6&  6&  6&  6&  6&  6&  6&  6&  6&  6& 6 \end{array} \right ] }\]

Now let's define the {\em $a$-loading} of a partition (this is closely related to PC1). By the Perron--Frobenius theorem, the eigenvector $v$ of $\mathsf Y_n$ corresponding to the largest eigenvalue will have all positive coordinates. As is common in statistics, we normalize the eigenvector as follows.
Denote its largest coordinate by $v_{\mathrm{max}}$ and the smallest by $v_{\mathrm{min}}$. Then, given a partition $\lambda$, we define its $a$-loading
as
\be
a_\lambda \seteq 100 \times \frac{v_\lambda - v_{\mathrm{min}}}{v_{\mathrm{max}} - v_{\mathrm{min}}} \quad \text{for } \lambda \in \mathcal P(n) .
\ee
For example, when $n$ is equal to 6, the $a$-loadings of partitions are given by the sequence
\begin{align*}
& (a_\lambda)_{\lambda \in \mathcal P(n)} = \\ & \phantom {LLL} (100.00,85.89,71.79,66.66,57.68,52.55,45.23,33.32,31.12,22.81,0.00).
\end{align*}
The partition $(6,0,0,0,0,0)$ has $a$-loading 100, while the partition $(1,1,1,1,1,1)$ has $a$-loading 0.

\subsection{Difference matrix and $b$-loadings}

Here is another natural embedding of the set of partitions into a higher dimensional space. 
We define a symmetric matrix $\mathsf Z_n$ whose $(\lambda, \mu)$-entry is given by the sum of the differences between the parts of the corresponding partitions:
\be
\mathsf {(Z_n)}_{\lambda, \mu}= \lVert \lambda-\mu \rVert_1 \seteq \sum_{i=1}^n | \lambda_i - \mu_i| 
\ee
for $\lambda=(\lambda_1, \lambda_2, \dots , \lambda_n)$ and $\mu=(\mu_1,\mu_2, \dots , \mu_n) \in \mathcal P(n)$.
For example,
\be
{\scriptsize \mathsf Z_6 = \left [ \begin{array}{ccccccccccc}  0  &2  &4  &4  &6  &6  &6  &8  &8  &8 &10 \\ 2  &0  &2  &2  &4  &4  &4  &6  &6  &6  &8 \\  4  &2  &0  &2  &2  &2  &4  &4  &4  &6  &8\\4  &2  &2  &0  &4  &2  &2  &4  &4  &4  &6\\6  &4  &2  &4  &0  &2  &4  &4  &4  &6  &8\\6  &4  &2  &2  &2  &0  &2  &2  &2  &4  &6\\6  &4  &4  &2  &4  &2  &0  &4  &2  &2  &4\\8  &6  &4  &4  &4  &2  &4  &0  &2  &4  &6\\8  &6  &4  &4  &4  &2  &2  &2  &0  &2  &4\\8  &6  &6  &4  &6  &4  &2  &4  &2  &0  &2\\ 10  &8  &8  &6  &8  &6  &4  &6  &4  &2  &0 \end{array} \right  ] }
\ee

Similarly to what we did above, we use PC1 of the matrix $\mathsf Z_n$ to define the {\em $b$-loadings} of a partition. Let $w$ be an eigenvector of the largest eigenvalue of $\mathsf  Z_n$ with all positive coordinates. Denote the largest coordinate by $w_{\mathrm{max}}$ and the smallest by $w_{\mathrm{min}}$. For a partition $\lambda$, define 
 \be 
 b_\lambda \seteq 100 \times \frac{w_\lambda - w_{\mathrm{min}}}{w_{\mathrm{max}} - w_{\mathrm{min}}} \quad \text{for } \lambda \in \mathcal P(n) .
 \ee
For example, when $n$ is equal to 6, the $b$-loadings of partitions are given by 
\begin{align*}
& (b_\lambda)_{\lambda \in \mathcal P(n)} = \\ & \phantom{LLL} (100.00,37.25,19.93,4.36,43.01,0.00,4.36,43.01,19.93,37.25, 100.00).
\end{align*}
The partition $(6,0,0,0,0,0)$ has $b$-loading 100, and so does the partition $(1,1,1,1,1,1)$. The $b$-loading 0 appears in the middle at $(3,2,1,0,0,0)$. 

\subsection{Loadings of triples}
 
To relate loadings to Kronecker coefficients, we consider triples of partitions $\mathbf t=(\lambda, \mu, \nu) \in \CP(n)^3$. 
For convenience, we will write $g(\mathbf t)$ to denote $g_{\lambda, \mu}^{\nu}$. The $a$-loading of a triple is defined to be the sum of the $a$-loadings of the individual partitions $a(\mathbf t) \seteq a_\lambda + a_\mu + a_\nu$.
Similarly, the $b$-loading of a triple is the sum of the $b$-loadings of each partition, $b(\mathbf t) \seteq b_\lambda + b_\mu + b_\nu$.

If we sketch the histograms of the $a$-loadings for $n=15,16$, we obtain these pictures. We conjecture that as $n\rightarrow\infty$ they are normally distributed, and already $n=15$ is very close to normal.

\begin{figure}[h!!!]
\begin{center}
\includegraphics[scale=0.29]{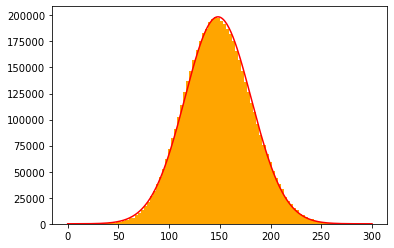}\quad \includegraphics[scale=0.29]{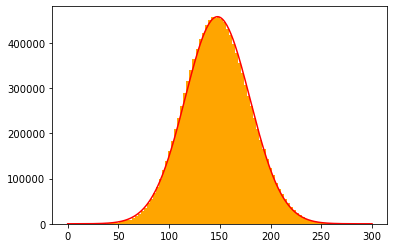}  
\end{center}
{\scriptsize \sf Histograms of $a$-loadings of $\mathbf t \in \mathcal P(n)^3$ for $n=15,16$} 
\end{figure}

On the other hand, the histograms of the $b$-loadings are these. We conjecture that its $n\rightarrow \infty$ limit is a gamma distribution, and again the fit for $n=15,16$ is very good.

\begin{figure}[h!!!]
\begin{center}
\includegraphics[scale=0.29]{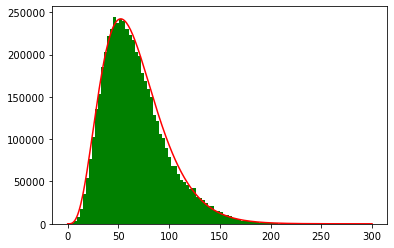}\quad \includegraphics[scale=0.29]{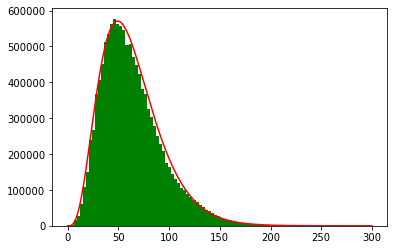}
\end{center}
{\scriptsize \sf Histograms of $b$-loadings of $\mathbf t \in \mathcal P(n)^3$ for $n=15,16$} 
\end{figure}

\subsection{Distributions conditioned on $g(\mathbf t) \neq 0$ or $g(\mathbf t) =0$} 

When we examine the distribution of $a$-loadings conditioned on whether the Kronecker coefficients are zero or non-zero, we get the following pictures. The distributions are different, but there is no clear separation. 

\begin{figure}[h!!!]
\begin{center}
 \includegraphics[scale=0.35]{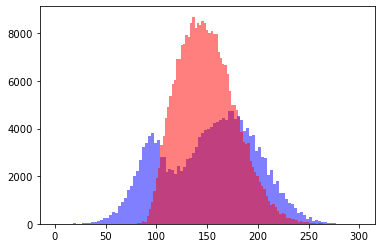} \quad \includegraphics[scale=0.35]{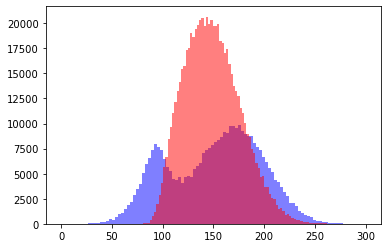}
\end{center}
{\scriptsize \sf Histograms of $a$-loadings of $\mathbf t \in \mathcal P(n)^3$ for $n=12,13$. The red (resp. blue) region represents the numbers of $\mathbf t$ such that $g(\mathbf t) \neq 0$ (resp. $g(\mathbf t) =0$).}
\end{figure}

On the other hand, the distribution of $b$-loadings exhibits separation, which can be used to confirm $g(\mathbf t) \neq 0$ for the separated regions.

\begin{figure}[h!!!]
\begin{center}
 \includegraphics[scale=0.35]{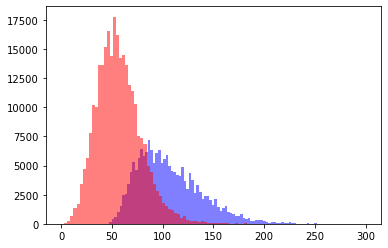}
\quad \includegraphics[scale=0.35]{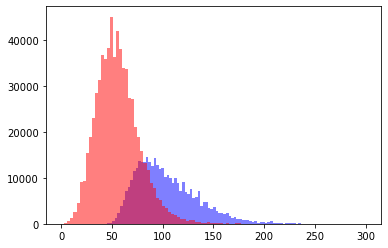} 
\end{center}
{\scriptsize \sf Histograms of $b$-loadings of $\mathbf t \in \mathcal P(n)^3$ for $n=12,13$. The red (resp. blue) region represents the numbers of $\mathbf t$ such that $g(\mathbf t) \neq 0$ (resp. $g(\mathbf t) =0$).}
\end{figure}

Define  
\begin{align*} 
b_{\star} &\seteq \min \{ b(\mathbf{t}) : g(\mathbf{t}) = 0 , \mathbf{t} \in \mathcal{P}(n)^3 \} .
\end{align*}  
Clearly, by definition, we have:  
\[
\text{For } \mathbf{t} \in \mathcal{P}(n)^3, \text{ if } b(\mathbf{t}) < b_{\star}, \text{ then } g(\mathbf{t}) \neq 0.
\]
While evident, {\it a priori} this could be vacuous, but as the graphs suggest it is not. 
For example, when $n=18$, consider \[ \mathbf t=((12, 4, 2), (8, 4, 2, 2, 1, 1), (5, 4, 3, 3, 1, 1, 1)).\]  Its $b$-loading is approximately $41.07$ $< b_{\star} \approx 44.18$.  Thus we immediately conclude that $g(\mathbf t) \neq 0$.  

How strong is it?
When $n=20$, there are 246,491,883 triples $\mathbf t$. Among them, 78,382,890 satisfy $b(\mathbf t) < b_\star$. This is about 31.8\%.
This is studied further in \cite{Kr}.  Although the picture is not totally clear, it supports
\begin{conjecture}
    The fraction of triples satisfying $b(\mathbf t) < b_\star$ has a nonzero infimum limit as $n\rightarrow\infty$.
\end{conjecture}
This is new as far as we know.
Some further studies are discussed in the paper of Butbaia, Lee and Ruehle in this volume:
the use of KANs (Kolmogorov--Arnold networks), symbolic regression and small neural networks.


\section{Concluding remarks}

This concludes our two case studies.   What can we take from them?  More specifically, what can we say about steps (3) and (4)
in \S \ref{ss:mds}, to interpret ML results to understand the mathematical objects better, and then
find new precisely defined structures, conjectures and theorems?  Here are a few observations (see also the many other works
in this volume).

First, while the primary goal of much scientific and practical ML is to model the data and make predictions, here that was only an initial step, essentially
to check whether the dataset actually contains any structure worth understanding.  The point is not so much to get 99\% or some other
accuracy figure in our fits, it is to check whether there are relations in the data which we can try to understand.  Many experiments fail this
initial test.  If so, one can try varying the choices of dataset and model discussed in \S \ref{ss:mds}.  If we are studying an infinite set of
objects, the choice of which finite subsets to use as datasets (the choice of size parameter or distribution to sample from) can be especially important, and 
negative results with one choice should not be overly discouraging.
But having made reasonable efforts, if these also fail one can move on.
This is the nature of exploration.

Second, in searching for structure in the data, the simplest methods are often the best.  In the case of the murmurations, simply assembling the data in the right way leads to the essential discovery.  As another example, while an FFN (feed forward neural network) might lead to a more accurate fit than a linear model, anything the linear model finds will be far easier to interpret.  Of course many problems cannot be treated by a linear model, for example to detect  an inequality which is nonlinear in the features.  But one should start simply.

Finally, like all mathematics to date, mathematical data science (MDS) is a human activity, bringing computers, statistics and data science to bear on diverse problems from many branches of mathematics, fundamentally very different but with aspects in common.  MDS is a way to identify and make use of some of these common aspects and to bring new concepts and tools to bear.  Perhaps in the future AI's will do mathematics in collaboration with humans or even on their own, and we predict that MDS (among many other mathematical approaches) will be one of the ways they do it.

\bibliographystyle{alpha}
\bibliography{bibfile}

\vskip 0.3 cm
\address{\hskip - 0.4 cm CMSA, Harvard University, Cambridge, MA 02138, USA}

\noindent \email{\href{mailto:mdouglas@cmsa.fas.harvard.edu}{mdouglas@cmsa.fas.harvard.edu}}
\vskip 0.3 cm

\address{\hskip - 0.4 cm Department of Mathematics, University of Connecticut, Storrs, CT 06269, USA  \hfill \break Korea Institute for Advanced Study, Seoul 02455, Republic of Korea}

\noindent \email{\href{mailto:khlee@math.uconn.edu}{khlee@math.uconn.edu}}

\end{document}